\documentclass{umi}
\usepackage{euscript}
\usepackage{amsfonts}
\usepackage{amssymb}
\usepackage{amsmath}
\usepackage{amscd}
\usepackage{xypic}
\usepackage{enumerate}
\newtheorem{theorem}{\sc Theorem}[section]
\newtheorem{lemma}[theorem]{\sc Lemma}

\newtheorem{remark}[theorem]{\sc Remark}

\newtheorem{proposition}[theorem]{\sc Proposition}
\newcommand{\C}{\mathbb C}

\newcommand{\PP}{\mathbb P}

\newcommand{\OO}{\mathcal{O}}

\newcommand{\GG}{\mathcal{G}}

\DeclareMathOperator{\HH}{H}
\DeclareMathOperator{\hh}{h}
\DeclareMathOperator{\Ext}{Ext}
\DeclareMathOperator{\Hom}{Hom}

\DeclareMathOperator{\rk}{rk}

\DeclareMathOperator{\GL}{GL}

\DeclareMathOperator{\End}{End}

\DeclareMathOperator{\I}{Id}
\DeclareMathOperator{\Stab}{Stab}
\begin{document}
%%%%%%%%%%%%%%%%%%%%%%%%%%%%%%%%%%%%%%%%%%%%%%%%%%%%%%%%%%%%%%
%
%               T  I  T  L  E
%
%%%%%%%%%%%%%%%%%%%%%%%%%%%%%%%%%%%%%%%%%%%%%%%%%%%%%%%%%%%%%
\title{Simplicity of generic Steiner bundles}
\author{{\sc Maria Chiara Brambilla}
\thanks{The author was supported by the MIUR in the framework of the National
 Research Projects `` Propriet\`a geometriche delle variet\`a reali e complesse'' and  ``Geometria delle
 variet\`a algebriche''.
}}
\date{}
\maketitle{}
\begin{abstract}
{Un fibrato di Steiner $E$ su $\PP^n$ ha una risoluzione lineare
della forma \linebreak
$0\rightarrow \OO(-1)^s\rightarrow \OO^t\rightarrow
E\rightarrow 0.$ In questo lavoro proviamo che il generico fibrato
di Steiner $E$ \`e semplice se e solo se $\chi(\End E)$ \`e minore
o uguale a $1$. In particolare mostriamo che  $E$ \`e eccezionale
oppure soddisfa la disuguaglianza $t\leq \left(
\frac{n+1+\sqrt{(n+1)^2-4}}{2}\right)s$.}
\end{abstract}
\begin{abstracting}
{A Steiner bundle $E$ on $\PP^n$ has a linear resolution of the
form \linebreak
$0\rightarrow \OO(-1)^s{\rightarrow} \OO^t\rightarrow
E\rightarrow 0.$ In this paper we prove that a generic Steiner
bundle $E$  is simple if and only if $\chi(\End E)$ is less or
equal to $1$. In particular we show that either $E$ is exceptional
or it satisfies the inequality $t\leq \left(
\frac{n+1+\sqrt{(n+1)^2-4}}{2}\right)s$.}
\end{abstracting}
%%%%%%%%%%%%%%%%%%%%%%%%%%%%%%%%%%%%%%%%%%%%%%%%%%%%%%%%%%%
%
%             I n t r o d u c t i o n
%
%%%%%%%%%%%%%%%%%%%%%%%%%%%%%%%%%%%%%%%%%%%%%%%%%%%%%%%%%%%%%%%
\section{Introduction}
According to \cite{DolgaKapra}
a Steiner bundle $E$ on $\PP(V)=\PP^{N-1}$
has a linear resolution of the form
$$0\rightarrow \OO(-1)^s{\rightarrow} \OO^t\rightarrow
E\rightarrow 0.$$ It is well known that Steiner bundles have rank
$t-s\geq N-1$ and if equality holds then they are stable, in
particular they are simple  (see \cite{AnconaOttaviani94}). The
aim of this paper is to investigate the simplicity of Steiner
bundles for higher rank.\\
\\
{\sc Main Theorem }
{\em Let $E$ be a Steiner bundle on $\PP^{N-1}$, with  $N\geq3$, defined by the exact sequence
$$0\rightarrow \OO(-1)^s\overset{m}{\rightarrow} \OO^t\rightarrow
E\rightarrow 0,$$
where $m$ is a generic morphism in $\Hom(\OO(-1)^s,\OO^t)$, then the
following statements are equivalent:
\begin{enumerate}
\item[\rm{(i)}]  $E$ is simple, i.e.\ $\hh^0(\End E)=1$,
\item[\rm{(ii)}] $s^2-Nst+t^2\leq 1$ i.e.\ $\chi(\End E)\leq 1$,
\item[\rm{(iii})] either  $s^2-Nst+t^2\leq 0$ i.e.\ $t\leq
(\frac{N+\sqrt{N^2-4}}{2})s$
  or
$(t,s)=(a_{k+1},a_k)$, where
$a_k=\frac{ \left(\frac{N+\sqrt{N^2-4}}{2}\right)^k-
  \left(\frac{N-\sqrt{N^2-4}}{2}\right)^k}{\sqrt{N^2-4}}.$
\end{enumerate}
} The generalized Fibonacci numbers appearing in $\rm{(iii)}$
satisfy a recurrence relation, as it is clear from the proof of
Theorem \ref{teoeccezionali}.

Our result in the case of $\PP^2$ is partially contained, although
somehow hidden, in \cite{DrezetLePotier}. Indeed Dr\'ezet and Le
Potier find a criterion to check the stability of a generic
bundle, given its rank and Chern classes. In the case of a
normalized Steiner bundle $E$ on $\PP^2$, it is possible to prove
that if $E$ satisfies condition $\rm{(iii)}$ of the main theorem,
then the Dr\'ezet-Le Potier condition for stability is satisfied.
Hence $E$ is stable and, consequently, simple. On the other hand,
when $E$ is not normalized, it is very complicated to check the
criterion of Dr\'ezet-Le Potier, but we can easily prove the
simplicity with other techniques. Anyway the proof that we present
in this paper is independent of \cite{DrezetLePotier}, is more
elementary and works on $\PP^n$ as well.

The genericity assumption cannot be dropped, because when
$\rk{E}=t-s>N-1$ it is always possible to find a decomposable
Steiner bundle, that is in particular non-simple.
%%%esempio 3,8,E+O

Since the equivalence between conditions $\rm{(ii)}$ and $
\rm{(iii)}$ is an arithmetic statement, our theorem claims that
$\chi(\End E)$ is the responsible for the simplicity of a generic
Steiner bundle $E$. Indeed it is easy to check that if $E$ is
simple then $\chi(\End E)\leq1$ (Lemma \ref{facile}) and this is
also true for some other bundles, for example for every bundle on
$\PP^2$. The converse is not true in general, because it is
possible to find a non-simple bundle $F$ on $\PP^2$ such that
$\chi(\End F)< 1$. For example we can consider the cokernel F of a
generic map of the form
$$0\rightarrow\OO(-2)\oplus\OO(-1)^4\rightarrow\OO^{16}\rightarrow
F\rightarrow 0,$$ where $\chi(\End F)=-3$, but it can be shown
that $\hh^0(\End F)=5$ therefore $F$ is not simple.

In the third statement of our theorem we claim that if $E$ is a simple
Steiner bundle,
 then either $E$ is
exceptional or it satisfies a numerical inequality (see Theorem
\ref{teoeccezionali}). We recall that exceptional bundles have no
deformations. The name exceptional in this setting is justified by
the fact that they are the only simple Steiner bundles which
violate the numerical inequality. It is remarkable to note that
all the exceptional bundles on $\PP^2$ can be constructed by the
theory of helices, in particular there exists a correspondence
between the exceptional bundles on the projective plane and the
solutions of the Markov equation $x^2+y^2+z^2=3xyz$ (see
\cite{Rudakov}).

The plan of the article is as follows:
 section $2$ is devoted to the case of exceptional bundles and section
$3$ to the proof of the main theorem.  At the end of the paper,
Theorem \ref{ananum} is a reformulation in terms of matrices of
the main theorem. As a basic reference for bundles on $\PP^n$ see
\cite{OSS}.

I would like to thank Giorgio Ottaviani, for suggesting me the
problem and for his continuous assistance, and Enrique Arrondo,
for many useful discussions. I also thank very much Jean Vall\`es,
for his helpful comments concerning this work, in particular for
his collaboration in simplifying the proof of Lemma
\ref{principale}. %%

%%%%%%%%%%%%%%%%%%%%%%%%%%%%%%%%%%%%%%%%%%%%%%%%%%%%%%%%%%%%%%%%%%%%%%%%%%%%%%%
%
%
%               Exceptional bundles
%
%
%%%%%%%%%%%%%%%%%%%%%%%%%%%%%%%%%%%%%%%%%%%%%%%%%%%%%%%%%%%%%%%%%%%%%%%%%%%%%%%
\section{Exceptional bundles}
In \cite{Rudakovseminari} the theory of
 helices of exceptional bundles is developed in  a general axiomatic presentation.
Here we give the following result as a particular case of this theory.

\begin{theorem}
\cite{Rudakov,Rudakovseminari}
\label{teoeccezionali}
Let $E_k$ be a generic Steiner bundle on $\PP^{N-1}$, with  $N\geq3$, defined by the exact sequence
$$0\rightarrow \OO(-1)^{a_{k-1}}{\rightarrow} \OO^{a_k}\rightarrow
E_k\rightarrow 0,$$
where
$$a_k=\frac{ \left(\frac{N+\sqrt{N^2-4}}{2}\right)^k-
  \left(\frac{N-\sqrt{N^2-4}}{2}\right)^k}{\sqrt{N^2-4}},$$
then $E_k$ is exceptional (i.e.\ $\hh^0(\End E)=1$ and $\hh^i(\End
E)=0$ for all $i>0$.)
\end{theorem}
%%%%%%%%%%%%%
%%%
On $\PP(V)=\PP^{N-1}$ we define a sequence of vector bundles as
follows:
\begin{equation}
\label{defF}
F_0=\OO(1), \quad F_1=\OO, \quad
F_{n+1}=\ker(F_n\otimes\Hom(F_n,F_{n-1})\overset{\psi_n}{\rightarrow}
F_{n-1}),
\end{equation}
where $\psi_n$ is the canonical map.
%
%%%%

The following lemma can be found in \cite{Rudakovseminari}. We
underline that it is possible to prove it in a straightforward way
only by standard cohomology sequences.
\begin{lemma}
\label{lemmaeliche} Given the definition \eqref{defF}, for all
$n\geq 1$  the canonical map $\psi_n$ is an epimorphism. Moreover
the following properties $(A_n)$, $(B_n)$ and $(C_n)$ are
satisfied for all $n\geq 1$:
\begin{eqnarray*}
(A_n)&\quad& \Hom(F_n,F_n)\cong\C,\quad \Ext^i(F_n,F_n)=0,\quad
\textrm{for all}\quad i\geq 1, \\
(B_n)&\quad& \Hom(F_{n-1},F_n)=0,\quad \Ext^i(F_{n-1},F_n)=0,\quad
\textrm{for all}\quad i\geq 1, \\
(C_n)&\quad &\Hom(F_n,F_{n-1})\cong V,\quad
\Ext^i(F_n,F_{n-1})=0,\quad\textrm{for all}\quad i\geq 1.
\end{eqnarray*}
Note that $(A_n)$ means that every $F_n$ is an exceptional bundle.
\end{lemma}
%%%%%%%%%%
%%%%
%%%%
\begin{remark}
Following  \cite{Rudakovseminari} the previous lemma means
 that $(F_n,F_{n-1})$ is a left admissible
pair and  $(F_{n+1},F_n)$ is the  left mutation %$L_{F_n}F_{n-1}$
of
$(F_n,F_{n-1})$  and that
the sequence $(F_n)$ forms an exceptional collection generated by the
helix $(\OO(i))$ by left mutations.
\end{remark}
%%%%%%%%%%%%%%%%%%%%%%%%%%%%%%%%%%%%%%%%%%%%%%%
{\sc Proof of Theorem \ref{teoeccezionali}.} Lemma
\ref{lemmaeliche} states that the bundles $F_n$, defined as in
\eqref{defF}, are exceptional for all $n\geq 0$. Obviously their
dual $F_n^*$ are exceptional too. Now we will prove that, for
every $n\geq1$, the bundle $F_n^*$ admits the following resolution
\begin{equation}
\label{claim}
0\rightarrow \OO(-1)^{a_{n-1}}{\rightarrow} \OO^{a_n}\rightarrow
F_n^*\rightarrow 0,
\end{equation}
 where $\{a_n\}$ is the sequence defined in the
statement. This implies that a generic bundle with this resolution
is exceptional. We can prove \eqref{claim} by induction on $n$.
First of all we notice that the sequence $\{a_n\}$  is also
defined recursively by
$$\left\{
\begin{array}{l}
a_0=0,\\
a_1=1,\\
a_{n+1}=Na_n-a_{n-1}.\\
\end{array}
\right.
$$
Therefore if $n=1$ the sequence \eqref{claim} is $0\rightarrow
\OO(-1)^{a_{0}}{\rightarrow} \OO^{a_1}\rightarrow F_1^*\rightarrow
0,$ i.e.\ $0{\rightarrow} \OO{\rightarrow} F_1^*\rightarrow 0,$
and this is true because  $F_1\cong \OO$. Now let us suppose that
every $F_k^*$ admits a resolution \eqref{claim} for all $k\leq n$
and we will prove it for $F_{n+1}^*$. Let us dualize the sequence
$$0\rightarrow F_{n+1}\rightarrow F_n\otimes\Hom(F_n,F_{n-1})\rightarrow
F_{n-1} \rightarrow 0
$$
and by induction hypothesis we have:
%%%%%%%%%%%%%%%%
%%%%%%%%%%%%%%
$$
\xymatrix{
       &      0          &  0                            &  &\\
0\ar[r]& F^*_{n-1}\ar[r]\ar[u]&   F^*_n\otimes V^*\ar[r]\ar[u]&   F^*_{n+1}\ar[r]&0 \\
       &\OO^{a_{n-1}}\ar[u]   &   \OO^{a_n}\otimes V^*  \ar[u]&      & \\
       &\OO(-1)^{a_{n-2}}\ar[u]&   \OO(-1)^{a_{n-1}}\otimes V^*  \ar[u]&  & \\
       &    0\ar[u]            &  0    \ar[u]            &  &
}
$$
We define the map $\alpha: \OO^{a_{n-1}}\rightarrow F^*_{n}\otimes
V^*$ as  the composition of the known maps. Since
$\Ext^1(\OO^{a_{n-1}},\OO(-1)^{a_{n-1}}\otimes
V^*)\cong\HH^1(\OO(-1))^{a_{n-1}^2}\otimes V^*=0$, the map
$\alpha$ induces a map $\widetilde\alpha:
\OO^{a_{n-1}}\rightarrow\OO^{a_n}\otimes V^*$
 such
that the following diagram commutes:
%%%%%%%%%%%%%%%%%%%%%%%%%%%%%%%%%%%%%%%%%
$$
\xymatrix{
       &        0        &  0                            &  &\\
0\ar[r]& F^*_{n-1}\ar[u]\ar[r]^f&   F^*_n\otimes V^*\ar[r]\ar[u]&
F^*_{n+1}\ar[r] &0 \\
       &        \OO^{a_{n-1}}\ar[u]   \ar@{-->}[r]^{\widetilde{\alpha}}\ar[ur]^{\alpha}       &   \OO^{a_n}\otimes V^*  \ar[u]&    & \\
       &          \OO(-1)^{a_{n-2}}\ar[u]        &   \OO(-1)^{a_{n-1}}\otimes V^*  \ar[u]&     & \\
       &              0\ar[u]  &  0    \ar[u]            &  &
}
$$
%%%%%%%%%%%%%%%%%%%%%%%%%%%%%%%%%%%%%%%%%%%
We observe that $\widetilde\alpha$ is injective if and only if
$\HH^0(\widetilde\alpha)$ is injective and, since
$\HH^0(\widetilde\alpha)=\HH^0(f)$, they are injective.
Obviously the cokernel of $\widetilde\alpha$ is
$\OO^{Na_n-a_{n-1}}=\OO^{a_{n+1}}$.
Let $\widetilde\beta$ be the restriction of $\widetilde\alpha$ to
$\OO(-1)^{a_{n-2}}$. Then we can check that
$\widetilde\beta$ is injective, its cokernel is
$\OO(-1)^{Na_{n-1}-a_{n-2}}=\OO(-1)^{a_{n}}$ and the following diagram commutes:
$$
\xymatrix{
       &         0       &  0                            &  0&\\
0\ar[r]& F^*_{n-1}\ar[u]\ar[r]&   F^*_n\otimes V^*\ar[r]\ar[u]&F^*_{n+1}\ar[r]\ar[u] &0 \\
0\ar[r]& \OO^{a_{n-1}}\ar[u]   \ar@{-->}[r]^{\widetilde{\alpha}}\ar[ur]^{\alpha}       &   \OO^{a_n}\otimes V^*\ar[r]  \ar[u]&   \OO^{a_{n+1}}  \ar[r]\ar[u]  &0 \\
0\ar[r]& \OO(-1)^{a_{n-2}}\ar[u]  \ar@{-->}[r]^{\widetilde{\beta}}      &   \OO(-1)^{a_{n-1}}\otimes V^*  \ar[u]\ar[r]& \OO(-1)^{a_n} \ar[u]     \ar[r]& 0\\
       &     0\ar[u]  &  0    \ar[u]            & 0\ar[u] &
}
$$
It follows that $F^*_{n+1}$ has the resolution
$0\rightarrow \OO(-1)^{a_{n}}{\rightarrow} \OO^{a_{n+1}}\rightarrow
F^*_{n+1}\rightarrow 0$ and
this completes the proof of our theorem.
\qed
%%%%%%%%%%%%%%%%%
%%%%%%%%%%%%%%%%%%%%%%%%%%%%%%%%%%%%%%%%%%%%%%%%%%%%%%%%%%%%%%%%%%%%%%%%%%%%%%%
%
%
%                 DIMOSTRAZIONE
%
%
%%%%%%%%%%%%%%%%%%%%%%%%%%%%%%%%%%%%%%%%%%%%%%%%%%%%%%%%%%%%%%%%%%%%%%%%%%%%%%%
\section{Proof of the main theorem}
Let $E$ be given by the exact sequence on $\PP^{N-1}=\PP(V)$
\begin{equation}
\label{successione}
0\longrightarrow I\otimes\OO(-1)\overset{m}{\longrightarrow}W\otimes \OO\longrightarrow
E\longrightarrow 0,
\end{equation}
where $V$, $I$ and $W$ are complex vector spaces of dimension $N\geq3$, $s$
and $t$ respectively
and $m$ is a generic morphism.
If we fix a basis in each of the vector spaces $I$ and $W$, the
morphism $m$ can be represented by a $t\times s$ matrix $M$ whose
entries are linear forms.
Let us consider the natural action of $\GL(I)\times\GL(W)$ on the
space $$H=\Hom(I\otimes \OO(-1),W\otimes \OO)\cong V\otimes
I^\vee\otimes W,$$ i.e.\ the action
$$H\times\GL(I)\times\GL(W)\rightarrow H$$
$$(M,A,B)\mapsto A^{-1}MB.$$
When the pair $(A,B)$ belongs to the stabilizer of $M$, it
 induces a
morphism $\phi:E\rightarrow E$, such that the following diagram
commutes:
\begin{equation}
\label{diagramma}
\xymatrix{
0\ar[r]&I\otimes\OO(-1)\ar[r]^M\ar[d]^A &W\otimes\OO\ar[r]\ar[d]^{B}&E\ar[r]\ar[d]^{\phi}&0 \\
0\ar[r]&I\otimes\OO(-1)\ar[r]^M         &W\otimes \OO\ar[r]
&E\ar[r]&0
}
\end{equation}
%%%%%%%%%%%%%%%%%%%%%%%%%%%%%%%%%%%%%%%%
%
%    Dimostrazione  che 1---->2
%
%%%%%%%%%%%%%%%%%%%%%%%%%%%%%%%%%%%%%%%%%%%%%%%%%%%%%%%%%%%%%%%%%%
%
\paragraph{I.}
Now we prove the first part of the theorem, i.e.\ the fact that
${\rm(i)}$ implies  ${\rm(ii)}$. %%
%%%%
\begin{remark}
From the sequence \eqref{successione} it follows that $\chi(E)=t$
and $\chi(E(1))=(Nt-s)$. Dualizing \eqref{successione} and
tensoring by $E$ we get
\begin{equation}
\label{succ-end1}
0\longrightarrow \End E\longrightarrow W^{\vee}\otimes E{\longrightarrow}
I^{\vee}\otimes E(1)\longrightarrow 0,
\end{equation}
therefore
$$\chi(\End E)=t\chi(E)-s\chi(E(1))=t^2-s(Nt-s)=t^2-Nst+s^2.$$
\end{remark}
%%%%%%%%
%%%%%%
%%%%%%
%%
\begin{lemma}
\label{facile}
If $E$ is a simple Steiner bundle, then $\chi(\End E)\leq1$.
\end{lemma}
\Proof From the sequences \eqref{successione} and
\eqref{succ-end1} it is easy to check that $\HH^i(\End E)=0$, for
all $i\geq2$. Moreover $\hh^0(\End E)=1$ because of the
simplicity, and consequently $\chi(\End E)=1-\hh^1(\End E)\leq1.$
\qed
%%%%%%%%%%%%%%%%%%%%%%%%%%%%%%%%%%%%%%%
%
%    Dimostrazione  che 2<---->3
%
%%%%%%%%%%%%%%%%%%%%%%%%%%%%%%%%%%%%%%%%%%%%%%%%%%%%%%%%%%%%%%%%%%
%
%
\paragraph{II.}
Now we prove that statement ${\rm(ii)}$ is equivalent to
${\rm(iii)}$.
\begin{remark}
Obviously  $s^2-Nst+t^2\leq 0$ is equivalent to
$(\frac{N-\sqrt{N^2-4}}{2})s \leq t\leq (\frac{N+\sqrt{N^2-4}}{2}) s$.
Since $t>s$ and $N>2$
%%%
this inequality is equivalent to $t\leq (\frac{N+\sqrt{N^2-4}}{2})s$.
Thus we have only to prove that  $s^2-Nst+t^2= 1$ is equivalent to
$(t,s)=(a_{k+1},a_k)$ where $a_k$ has been defined above.
\end{remark}
%%%
%%%
\begin{lemma}
All the integer solutions of $s^2-Nst+t^2=1$, when $t>s$, are
exactly $s=a_k,t=a_{k+1}$, where $a_k=\frac{
\left(\frac{N+\sqrt{N^2-4}}{2}\right)^k-
  \left(\frac{N-\sqrt{N^2-4}}{2}\right)^k}{\sqrt{N^2-4}}.$
\end{lemma}
\Proof We already know that the sequence $\{a_k\}$ is defined
recursively by
$$\left\{
\begin{array}{l}
a_0=0,\\
a_1=1,\\
a_{k+1}=Na_k-a_{k-1}.\\
\end{array}
\right.
$$
So we prove by induction on $k$ that $(s=a_k,t=a_{k+1})$ is a
solution of
\begin{equation}
\label{equazione}
s^2-Nst+t^2=1.
\end{equation}
If $k=0$, obviously $(s=0,t=1)$ is a solution. Let the pair
$(a_{k-1},a_k)$ satisfy \eqref{equazione}, then, using the
recursive definition, we check that $(a_{k},a_{k+1})$ is a
solution too. Hence we have to prove that there are no other
solution. By the change of coordinates $\{r=2t-Ns, s=s \}$ our
equation becomes the following Pell-Fermat equation
$r^2-(N^2-4)s^2=4$. By Number Theory results (see for example
\cite{Samuel}, page $77$, or \cite{Lenstra}), we know that all the
solutions $(r,s)$ are given by the sequence $(r_k,s_k)$ defined by
$$r_k+s_k\sqrt{N^2-4}=\frac{1}{2^{k-1}}(N+\sqrt{N^2-4})^k,$$
for all $k\geq0$. Now we have only to prove that these solutions
are exactly those already known. We can easily check that the pair
of sequences $(s_k,t_k)$ can be recursively defined by
$$\left\{
\begin{array}{l}
r_0=2,\\
s_0=0,\\
r_{k+1}=\frac{(N^2-4)s_k+Nr_{k}}{2},\\
s_{k+1}=\frac{Ns_k+r_k}{2}.\\
\end{array}
\right.
$$
By a change of coordinates we define $t_k=\frac{Ns_k+r_k}{2}$ and
we check that the pair $(s_k,t_k)$ is exactly $(a_k,a_{k+1})$, for
all $k\geq0$. In fact $(s_0,t_0)=(0,1)=(a_0,a_1)$ and, moreover,
$t_k=s_{k+1}$ and
$t_{k+1}=\frac{Ns_{k+1}+r_{k+1}}{2}=\frac{{(N^2-2)s_k+Nr_{k}}}{2}=
\frac{{(N^2-2)s_k+N(2t_k-Ns_k)}}{2}=Nt_k-t_{k-1}.$ \qed
%%%%%%%%%%%%%%%%%%%%%%%%%%%%%%%%%%%%%%%%%%%%%%%%%%%%%%%%%%%%%%%%%%%%%%%%%%%%%%%%%
%
%            Dimostrazione che 3--->1
%
%%%%%%%%%%%%%%%%%%%%%%%%%%%%%%%%%%%%%%%%%%%%%%%%%%%%%%%%%%%%%%%%%%%
\paragraph{III.}
Now we prove  the last implication, i.e.\ $\rm{(iii)}$ implies
$\rm{(i)}$. In the case $(t,s)=(a_{k+1},a_k)$, the generic $ E$ is
an exceptional bundle by Theorem  \ref{teoeccezionali}, therefore
it is in particular simple. So suppose $s^2-Nst+t^2\leq 0$ and
%%%%%%%%
recall that $H$ denotes $\Hom(I\otimes \OO(-1),W\otimes\OO)\cong
V\otimes I^\vee\otimes W.$
Let $S$  be the set
$$\{A,B,M : A^{-1}MB=M\}\subset \GL(I)\times\GL(W)\times H$$
and $\pi_1$ and $\pi_2$ the projections on $ \GL(I)\times\GL(W)$
and on $H$ respectively. Notice that, for all $M\in H$,
$\pi_1(\pi_2^{-1}(M))$ is  the stabilizer of $M$ with respect to
the action of $ \GL(I)\times\GL(W)$. Obviously $(\lambda\I
,\lambda\I)\in\Stab(M)$, therefore $\dim\Stab(M)\geq1$.
%%%%%%%
\begin{lemma}
\label{lemmastab1}
 If $E$ is defined by the sequence
\begin{equation}
\label{equazM}
0\longrightarrow I\otimes\OO(-1)\overset{M}{\longrightarrow}W\otimes
\OO\longrightarrow E\longrightarrow 0
\end{equation}
and $\dim \Stab(M)=1$, then $E$ is simple.
\end{lemma}
\Proof If by contradiction  $E$ is not simple, then there exists
$\phi:E\rightarrow E$ non-trivial. Applying the functor
$\Hom(-,E)$ to the sequence \eqref{equazM}
we get that $\phi$ induces $\widetilde{\phi}$ non-trivial in
$\Hom(W\otimes \OO,E)$. Now applying the functor
$\Hom(W\otimes\OO,-)$ again to the same sequence  we get
 $\Hom(W\otimes\OO,W\otimes\OO)\cong \Hom(W\otimes\OO,E)$
because $\Hom(W\otimes\OO, I\otimes\OO(-1))\cong W\otimes I\otimes
\HH^0(\OO(-1))=0$ and $\Ext^1(W\otimes\OO,I\otimes\OO(-1))\cong
W\otimes I\otimes \HH^1(\OO(-1))=0$. It follows that there exists
$\widetilde{\phi}$ non-trivial in $\End(W\otimes\OO)$, i.e.\ a
matrix $B\neq\I$ in $\GL(W)$. Restricting  $\widetilde{\phi}$ to
$I\otimes\OO(-1)$ and calling $A$ the corresponding matrix in
$\GL(I)$, we get the commutative diagram \eqref{diagramma}.
Therefore $(A,B)\neq(\lambda\I,\lambda\I)$ belongs to $\Stab(M)$
and consequently $\dim\Stab(M)>1$. \qed
%%%%%%%%%%%
Finally it suffices to prove that for all generic $M\in H$, the dimension
of the stabilizer is exactly $1$.
In other words we have to prove the following
\begin{proposition}
\label{proposizione}
Let $H=V\otimes I^\vee\otimes W$ as above and suppose $s^2-Nst+t^2\leq 0$.
Then the generic orbit in
$H$ with respect to the natural action of
$GL(I)\times\GL(W)$ has dimension exactly $(s^2+t^2-1)$.
\end{proposition}
%%%%%%
Recall that we have defined the following diagram
$$
\xymatrix@C-10ex{
  & S=  \{A,B,M : A^{-1}MB=M\}\ar[dl]^{\pi_1}\ar[dr]_{\pi_2}&\\
\GL(I)\times\GL(W)&&H
}
$$
Let $(A,B)$ be two fixed Jordan canonical forms in $\GL(I)\times\GL(W)$.
We define $G_{AB}\subset\GL(I)\times\GL(W)$ as the set of couples of
matrices similar respectively to
$A$ and $B$.
Note that
$\pi_2\pi_1^{-1}(G_{AB})=\{C^{-1}MD:A^{-1}MB=M, C\in\GL(I), D\in\GL(W)\}$.
Moreover
$G_{\I\I}=\{(\lambda\I,\lambda\I),\lambda\in\C\}$
and $\pi_2\pi_1^{-1}(G_{\I\I})=H$.
\begin{lemma}\label{principale}
If $s^2-Nst+t^2\leq 0$ and $(A,B)$ are Jordan canonical forms
different from  $(\lambda\I,\lambda\I)$ for any $\lambda$, then
$\pi_2\pi_1^{-1}(G_{AB})$ is contained in a Zariski closed subset
strictly contained in $H$.
\end{lemma}
\Proof Suppose that the assertion is false. Then there exist two
Jordan canonical forms $A$ and $B$, different from
$(\lambda\I,\lambda\I)$, such that $\pi_2\pi_1^{-1}(G_{AB})$ is
not contained in any closed subset. This implies that we can take
a general $M\in H$ such that $AM=MB$ and in particular we can
suppose the rank of $M$ maximum.

Now we prove that $A$ and $B$ have the same minimal polynomial.
First, if $p_B$ is the minimal polynomial of $B$, i.e.\
$p_B(B)=0$, then it follows that $p_B(A)M=Mp_B(B)=0$ and since $M$
is injective we get $p_B(A)=0$, hence  the minimal polynomial of
$B$ divides that of $A$. Now if we denote by $\lambda_i$ ($1\leq
i\leq q$) the eigenvalues of $A$ and by $\mu_j$ ($1\leq j\leq q'$)
those of $B$, we obtain that
$\mu_j\in\{\lambda_1,\ldots,\lambda_q\}$ for all $1\leq j \leq
q'$. Let us define $A'=(A-x\I_s)$ and $B'=(B-x\I_t)$: obviously we
obtain $A'M=MB'$. We denote by $\overline{B'}$ the matrix of
cofactors of $B'$ and we know that
$B'\overline{B'}=\det(B')\I_t=P_B(x)\I_t$, where $P_B$ is the
characteristic polynomial of $B$. Therefore
$$A'M\overline{B'}=P_B(x)M$$ and developing this expression we see
that $q'=q$. In fact if there exists a $\lambda_i\neq\mu_j$ for
all $j=1,\ldots,q'$, then there is  a row of zeroes in $M$ and
consequently $M$ is not generic.
%%%%%%
%%%%%%%%%%%%%%
Then we get $A$ and $B$ with the same eigenvalues $\lambda_i$
($1\leq i\leq q$) with multiplicity respectively $a_i\geq 1$  and
$b_i\geq 1$. The hypothesis that $(A,B)\neq(\lambda\I,\lambda\I)$
means that either $A$ and $B$ have more than one eigenvalue or at
least one of them is non-diagonal.

Now consider the first case, i.e.\ $q\geq2$. Since $\dim I= s$ and
$\dim W= t$, obviously $\sum_{i=1}^q a_i=s$ e $\sum_{i=1}^q
b_i=t$. Now we denote  $M=(M_{ij})$, where $M_{ij}$ has dimension
$a_i\times b_j$. Since $AM=MB$, every block  $M_{ij}$ is zero for
all $i\neq j$, i.e.\ it is possible to write $M$ with the form
$$   M=
\left(\begin{array}{cccc}
* &0&\cdots&0 \\
0&*&\cdots&0\\
0&0&\ddots&0\\
0&0&\cdots&*
\end{array}
\right).$$
In particular we can define $n_1=a_1, n_2=\sum_{i=2}^p a_i, m_1=b_1, m_2=\sum_{i=2}^p b_i$
and thus
the matrix $M$ becomes
\begin{equation}
\label{matrice}
M=
\left(\begin{array}{cc}
(*)_{n_1\times m_1}&(0)_{n_1\times m_2}\\
(0)_{n_2\times m_1}&(*)_{n_2\times m_2}
\end{array}
 \right)
\end{equation}
where $n_1+n_2=s$ and $m_1+m_2=t$ and $n_i, m_i\geq1$ for $i=1,2$.
Thus it only suffices to show that a matrix in the orbit
$$O_M=\{C^{-1}MD: C\in\GL(s), D\in\GL(t), M \textrm{ with the form $\eqref{matrice}$}\},$$
is not generic in $H$ if $s^2-Nst+t^2\leq0$. This fact contradicts
our assumption and completes the proof.

In order to show this, we introduce the following diagrams
$$
\xymatrix@C-10ex{
  &   \{\phi,I_1,W_1:\phi(I_1\otimes V^\vee)\subseteq W_1\}\ar[dl]^{\alpha_1}\ar[dr]_{\beta_1}&\\
H=\Hom( I \otimes V^\vee , W )&& \GG_1=\GG(\C^{n_1},\C^{s})\times\GG(\C^{m_1},\C^{t}))}
$$
where $\GG(\C^k,\C^h)$
denotes the Grassmannian of $\C^k\subset\C^h$ and
$$
\xymatrix@C-10ex{
  &   \{\phi,I_2,W_2:\phi(I_2\otimes V^\vee)\subseteq W_2\}\ar[dl]^{\alpha_2}\ar[dr]_{\beta_2}&\\
H=\Hom( I \otimes V^\vee , W )&& \GG_2=\GG(\C^{n_2},\C^{s})\times\GG(\C^{m_2},\C^{t}))}
$$
It is easy to check that the matrices of the
set $O_M$ live in the subvariety
$$\widetilde{H}=\alpha_1(\beta_1^{-1}(\GG_1))\cap\alpha_2(\beta_2^{-1}(\GG_2))\subseteq
H,$$
then, in order to prove that these matrices are not generic, it
suffices to show that $\dim \widetilde{H}<\dim H$.
Since $\dim(\GG_i)=(n_1n_2+m_1m_2)$ for $i=1,2$,
we obtain
$$\dim(\alpha_1(\beta_1^{-1}(\GG_1)))\leq\dim(\beta_1^{-1}(\GG_1))=n_1n_2+m_1m_2+N(n_1(m_1+m_2)+n_2m_2)$$
and
 $$\dim(\alpha_2(\beta_2^{-1}(\GG_2)))\leq\dim(\beta_2^{-1}(\GG_2))=n_1n_2+m_1m_2+N(n_1m_1+n_2(m_1+m_2)).$$
Therefore, since $\dim H=Nst=N(n_1+n_2)(m_1+m_2)$
 we only need to show that either
$(n_1n_2+m_1m_2-Nn_2m_1)<0$ or
$(n_1n_2+m_1m_2-Nn_1m_2)<0$.
In other words we have to prove that the system
$$\left\{
\begin{array}{l}
n_1n_2+m_1m_2-Nn_1m_2\geq 0
\\
n_1n_2+m_1m_2-Nn_2m_1\geq 0
\end{array}\right.
$$
has no solutions in our hypotesis $s^2-Nst+t^2\leq 0$, i.e.\ if
$$\frac{N-\sqrt{N^2-4}}{2}t\leq s\leq\frac{N+\sqrt{N^2-4}}{2}t.$$
This is equivalent to prove that the system
$$
\left\{
\begin{array}{l}
n_1n_2+m_1m_2-Nn_1m_2\geq 0
\\
n_1n_2+m_1m_2-Nn_2m_1\geq 0
\\
n_1+n_2\geq\frac{N-\sqrt{N^2-4}}{2}(m_1+m_2)
\\
n_1+n_2\leq\frac{N+\sqrt{N^2-4}}{2}(m_1+m_2)
\end{array}\right.
$$
has no solutions.
In order to do it,
consider  $n_1$ and $m_1$ as parameters and
write the previous system as a system of linear
inequalities in two unknowns $n_2$ and $m_2$:
$$
\left\{
\begin{array}{l}
n_1n_2\geq(Nn_1-m_1)m_2
\\
(n_1-Nm_1)n_2\geq-m_1m_2
\\
n_2\geq\alpha_-m_2+(\alpha_-m_1-n_1)
\\
n_2\leq\alpha_+m_2+(\alpha_+m_1-n_1)
\end{array}\right.
$$
where we denote $\alpha_-=\frac{N-\sqrt{N^2-4}}{2}$ and
$\alpha_+=\frac{N+\sqrt{N^2-4}}{2}$. Notice  that
$(\alpha_-+\alpha_+)=N$ and $\alpha_-\alpha_+=1$, because they are
solutions of the equation $s^2-Nst+t^2=0$. Now let us consider
three cases:
\begin{itemize}
\item if $0<n_1-\alpha_+m_1$
the system
$$
\left\{
\begin{array}{l}
n_2\geq\frac{(Nn_1-m_1)}{n_1}m_2
\\
n_2\leq\alpha_+m_2+(\alpha_+m_1-n_1)
\end{array}\right.
$$
has no solutions because $(\alpha_+m_1-n_1)<0$ and
$\alpha_+<\frac{(Nn_1-m_1)}{n_1}$, since
$(N-\alpha_+)n_1-m_1=\alpha_-n_1-m_1=(\alpha_+)^{-1}(n_1-\alpha_+m_1)>0$;
\item if $ n_1-\alpha_+m_1<0<n_1-\alpha_-m_1$  the system is
$$
\left\{
\begin{array}{l}
n_2\geq\frac{(Nn_1-m_1)}{n_1}m_2
\\
n_2\leq\frac{m_1}{(Nm_1-n_1)}m_2
\end{array}\right.
$$
because $Nm_1-n_1>\alpha_+m_1-n_1>0$ and there is no
solution because
$\frac{m_1}{(Nm_1-n_1)}<\frac{(Nn_1-m_1)}{n_1}$, since
$N(Nn_1m_1-m_1^2-n_1^2)>0$;
\item if $n_1-\alpha_-m_1<0$ then the system
$$
\left\{
\begin{array}{l}
n_2\leq\frac{m_1}{(Nm_1-n_1)}m_2
\\
n_2\geq\alpha_-m_2+(\alpha_-m_1-n_1)
\end{array}\right.
$$
has no solutions because $(\alpha_-m_1-n_1)>0$ and
$\alpha_->\frac{m_1}{(Nm_1-n_1)}$ i.e.\
$\alpha_+<\frac{(Nm_1-n_1)}{m_1}$, since
$(N-\alpha_+)m_1-n_1=\alpha_-m_1-n_1>0$.
\end{itemize}
Thus the  proof in the case $q\geq2$ is complete.

In the second case we consider $q=1$ and the two matrices are
$$   A=
\left(\begin{array}{ccc}
J_1& &\\
&\ddots&\\
&&J_h
\end{array}\right),
\quad
{\textrm{ where }}\quad
 {J_i=
\left(\begin{array}{cccc}
\lambda&1 & &\\
&\lambda&1&\\
&&\ddots&\ddots\\
&&&\lambda\\
\end{array}\right)}$$
and  $c_i$ denotes the order of $J_i$
and
$$   B=
\left(\begin{array}{ccc}
L_1& &\\
&\ddots&\\
&&L_k
\end{array}\right),
\quad
{\textrm{ where }}\quad
 {L_i=
\left(\begin{array}{cccc}
\lambda&1 & &\\
&\lambda&1&\\
&&\ddots&\ddots\\
&&&\lambda\\
\end{array}\right)}$$
and $d_i$ is the order of $L_i$. We suppose that $c_1\geq2$ or
$d_1\geq2$ i.e.\ $h<s$ or $k<t$. Then a matrix $M$ such that
$AM=MB$ has the form $  M=(M_{ij})$ and $M_{ij}$ is a $c_i\times
d_j$ matrix such that
$$M_{ij}=\left\{
\begin{array}{ccc}
T_c&{\textrm{ if }}&c_{i}=d_{j}=c\\
(0| T_c)&{\textrm{ if }}&c=c_{i}<d_{j}\\
(\frac{T_d}{0})&{\textrm{ if }}&c_{i}>d_{j}=d
\end{array}\right.
$$
and $T_c$ is a $c\times c$  upper-triangular Toeplitz matrix. It
is easy to see that $M$ has at least $k$ columns in which there
are at least $(c_1-1)+(c_2-1)+\ldots+(c_h-1)=(s-h)$ zeroes in such
a way that we can order the basis so as to write $M$ in the
following form
$$
\left(\begin{array}{cc}
(*)_{h\times k}&(*)_{h\times (t-k)}\\
(0)_{(s-h)\times k}&(*)_{(s-h)\times (t-k)}
\end{array}
 \right).
$$
Analogously $M$ has at least $h$ rows with at least $(t-k)$ zeroes such that it
is possible to write the matrix in the form
$$
\left(\begin{array}{cc}
(*)_{h\times k}&(0)_{h\times (t-k)}\\
(*)_{(s-h)\times k}&(*)_{(s-h)\times (t-k)}
\end{array}
 \right).
$$
Hence there exist non-trivial subspaces $I_1,I_2,W_1,W_2$ such
that $M(I_i\otimes V^{\vee})\subseteq W_i$, for $i=1,2$, and $\dim
I_1=s-h$, $\dim W_1=k$, $\dim I_2=h$, $\dim W_2=t-k$. Therefore
exactly the same argument used in the first case gives that $M$ is
not generic and completes the proof. \qed The previous lemma
proves  Proposition \ref{proposizione} and the main theorem
follows. This theorem can also be reformulated as follows:
\begin{theorem}
\label{ananum}
Let $M$ a  $(s\times t)$ matrix whose
entries are linear forms in
$N$ variables and consider the system
\begin{equation}
\label{xmmy}
XM=MY,
\end{equation}
where $X\in\GL(s)$
and $Y\in\GL(t)$ are the unknowns.
Then if $s^2+t^2-Nst\leq 1$,
there is a dense subset of the vector space
$\C^s\otimes\C^t\otimes\C^N$, where $M$ lives, such that
 the only solutions
of \eqref{xmmy}
 are trivial, i.e.\
$(X,Y)=(\lambda\I,\lambda\I)\in\GL(s)\times\GL(t)$ for
$\lambda\in\C$. Conversely if $s^2+t^2-Nst\geq 2$, then for all
$M$ there are non-trivial solutions.
\end{theorem}
\vspace{1.5cm}
%%%%%%%%%%%%%%%%%%%%%%%%%%%%%%%%%%%%%%%%%%%%%%%%%%%%%%%%%%%%%%%%%%%%%%%%
%
%                 R  E  F  E  R  E  N  C  E  S
%
%%%%%%%%%%%%%%%%%%%%%%%%%%%%%%%%%%%%%%%%%%%%%%%%%%%%%%%%%%%%%%%%%%%%%%%%

%%%%%%%%%%%%%%%%%%%%%%%%%
\vspace{1cm}

\center{Dipartimento di Matematica ``U.Dini''\\
Viale Morgagni 67/A\\
50134 Firenze, Italy \\
email: brambilla@math.unifi.it}

\end{document}